\documentclass[oneside, 12pt]{amsart} \topmargin      -10mm \textwidth      160 true mm
\usepackage{amsmath, amssymb, amsfonts, amstext, amsthm, amscd}
\usepackage[mathscr]{euscript}
\usepackage{enumerate}
\usepackage{tikz,color,soul}
\usetikzlibrary{matrix,arrows}

\usepackage[utf8x]{inputenc}
\usepackage[T1]{fontenc}

\usepackage[english]{babel}

\usetikzlibrary{arrows.meta, positioning,arrows}
\newtheorem{thm}{Theorem}[subsection]
\newtheorem{prop}[thm]{Proposition}

\newtheorem{cor}[thm]{Corollary}
\newtheorem{defn}[thm]{Definition}
\newtheorem{ex}[thm]{Example}
\newtheorem{rmk}[thm]{Remark}

\numberwithin{equation}{section}

\newcommand{\C}{\mathbb C}

\newcommand{\N}{\mathbb N}

\newcommand{\R}{\mathbb R}
\newcommand{\Z}{\mathbb Z}

\begin{document}

\title{Toric vector bundles on Bott tower}

\author{Bivas Khan}
\address{Department of Mathematics, Indian Institute of Technology-Madras, Chennai, India}
\email{bivaskhan10@gmail.com}

\author{Jyoti Dasgupta}
\address{Department of Mathematics, Indian Institute of Technology-Madras, Chennai, India}
\email{jdasgupta.maths@gmail.com}



\begin{abstract}
In this paper, using Klyachko's classification theorem we study positivity and semi-stability of toric vector bundles on a class of nonsingular projective toric varieties, known as Bott towers. In particular, we give a criterion of $s$-jet ampleness of line bundles and characterize nef and big line bundles on Bott towers using Bott numbers. We obtain a criterion for the ampleness of discriminant zero semi-stable toric vector bundles on nonsingular projective toric varieties. We also describe toric subbundles of toric vector bundles on nonsingular toric varieties.  
\end{abstract}

\maketitle

\section{Introduction}
Let $X$ be a complex toric variety under the action of a complex algebraic torus $T$. A vector bundle $\pi : \mathcal{E} \rightarrow X$ is said to be a toric vector bundle on $X$ if it has an equivariant $T$-structure, i.e. an action of the torus on $\mathcal{E}$ which is linear on fibres and the morphism $\pi$ is $T$-equivariant. Toric vector bundles were first classified by Kaneyama for nonsingular complete toric varieties (see \cite{kane}). Later Klyachko gave a combinatorial description for arbitrary toric varieties. He showed that toric vector bundles on toric varieties are equivalent to a family of filtrations on a vector space satisfying certain compatibility conditions (\cite[Theorem 2.2.1]{kly}). Tangent bundles of nonsingular toric varieties arise as natural examples of toric vector bundles. Using Klyachko's description, we give a necessary and sufficient condition for equivariant splitting of the tangent bundle of nonsingular toric varieties (see Proposition \ref{18}). We give a complete description of toric subbundles of a given toric vector bundle (see Proposition \ref{6}). As an application, we study stability of tangent bundle on Hirzebruch surface which recovers the result obtained in \cite[Theorem 6.2]{ozhan} (see Corollary \ref{14}). 

We consider a class of smooth projective toric varieties, known as Bott towers. A Bott tower of height $n$ $$M_n \rightarrow M_{n-1} \rightarrow \cdots \rightarrow M_2 \rightarrow M_1 \rightarrow M_0=\{ \text{point} \}$$ is defined inductively as an iterated projective bundle so that at the $k$-th stage of the tower,  $M_k$ is of the form $\mathbb{P}(\mathcal{O}_{M_{k-1}}  \oplus \mathcal{L}$) for an arbitrarily chosen line bundle $\mathcal{L}$ over $M_{k-1}$. Bott towers were shown to be deformation of Bott-Samelson varieties by Grossberg and Karshon in \cite{grossbergkarshon}. A Bott tower of height $n$ is completely determined by $\frac{n(n-1)}{2}$ integers $\{c_{i, j}\}_{\{1\leq i < j \leq n \}}$, known as the Bott numbers (see subsection 2.2). Without loss of generality, we can assume that the Bott numbers are non-negative (see Theorem \ref{12}). In this paper we are interested in studying the positivity properties of toric vector bundles on Bott tower.

We give a criterion for $s$-jet ampleness of line bundles on Bott towers (see Theorem \ref{16}). As a corollary we deduce a criterion for ampleness and nefness of line bundles on Bott towers\footnote{While working on this project, it came to our notice that Narasimha B. Chary has obtained similar criterion for ample and nef line bundles on Bott towers in \cite{bonala2017mori} using a different method.} (see Corollary \ref{2}). We give a criterion for an invariant divisor on a nonsingular projective toric variety to be nef and big (see Proposition \ref{3}). Using this, we give a characterization of nef and big line bundles on Bott towers using Bott numbers (see Theorem \ref{5}), and consequently Fano (respectively, weak Fano) Bott towers are determined. 

In general, the ampleness of a vector bundle $\mathcal{E}$ cannot be read off from ampleness of its determinant bundle $\text{det}(\mathcal{E})$. For example, consider the Hirzebruch surface \( \mathcal{H}_1\) which is Fano (see Corollary \ref{2}), but its tangent bundle \(T_X\) is not ample by Mori's theorem (see \cite[Remark \(6.3.2\)]{lazarsfeld14500positivity}). However, to study ampleness of discriminant zero semi-stable toric vector bundles on nonsingular projective toric variety, it is enough to consider its determinant bundle (see Proposition \ref{7}). The proof makes use of the results of Hering, Musta{\c{t}}{\u{a}}, Payne \cite{hering2010positivity} and Bruzzo and Giudice \cite{bruzzo2013restricting}.

\noindent
{\bf Acknowledgements:} We thank our advisors Arijit Dey and V. Uma for their valuable guidance throughout this project. We also thank D.S. Nagaraj for his valuable comments and suggestions on an earlier version of the manuscript. We thank the Council of Scientific and Industrial Research (CSIR) for their financial support. We thank the anonymous referee for his/her insightful comments and suggestions which led to many improvements of the manuscript.
\section{Some basic results}

\subsection{Toric vector bundle}

Let $T \cong \left(\C^*\right) ^n$ be an algebraic torus. Let $M=\text{Hom}(T, \C^*) \cong \Z^n$ be its character lattice and $N=\text{Hom}_{\Z}(M, \Z)$ be the dual lattice. Let $\Delta$ be a fan in $N_{\R}:=N \otimes_{\Z} \R$ which defines a toric variety $X=X(\Delta)$ under the action of the torus $T$. Let $\Delta(1)$ denote the edges of $\Delta$ and $\sigma(1)$ denote the edges of a cone $\sigma$ in $\Delta$. Let \(x_{\sigma}\) denote the distinguished point of the affine toric variety \(U_{\sigma}\) corresponding to the cone \(\sigma\) (see \cite[Section 3.2, Page 116]{Cox}, \cite[Section 2.1]{Ful}). Let \(T_{\sigma}\) denote the stabilizer of \(x_{\sigma}\) and $\widehat{T}_{\sigma}$ denote its character lattice. For each $\rho \in \Delta(1)$, let $D_{\rho}$ denote the $T$-invariant prime divisor corresponding to $\rho$. \\
A vector bundle $\pi : \mathcal{E} \rightarrow X$ is a toric or $T$-equivariant vector bundle on $X$ if $\mathcal{E}$ has a lift of the action of $T$ such that:
\begin{enumerate}
	\item The action of $T$ on $\mathcal{E}$ is linear on fibers of $\pi$. 
	\item $\pi$ is $T$-equivariant i.e., for all $ e \in \mathcal{E}$ and $t \in T$, $\pi(t \cdot e) = t \cdot \pi(e)$. 
\end{enumerate}

\noindent
It is known that any line bundle on toric variety is isomorphic to a toric line bundle. In fact, a choice of $T$-invariant Cartier divisor $D$ such that $\mathcal{L} \cong \mathcal{O}_X(D)$ gives rise to an equivariant structure on a line bundle $\mathcal{L}$ on \(X\) (see \cite[Page 13]{jgonza}). 

\noindent
 Klyachko gave a classification in terms of filtrations of certain vector space as follows:

\begin{thm}\label{1} \cite[Theorem $2.2.1$]{kly} The category of toric vector bundles over the toric variety $X=X(\Delta)$ is equivalent to the category of vector spaces $E$, with a family of decreasing $\Z$-filtrations $E^{\rho}(i)$, for $\rho \in \Delta(1)$, which satisfy the following compatibility condition:\\
$({\bf C})$ for any $\sigma \in \Delta$, there exists a $\widehat{T}_{\sigma}$-grading $E=\oplus_{\chi \in \widehat{T}_{\sigma}}E^{[\sigma]} ({\chi})$, for which 
$$E^{\rho}(i)=\bigoplus_{\langle \chi, v_{\rho} \rangle \geqslant i}E^{[\sigma]}({\chi})$$ for all $\rho \in \sigma(1)$, where $v_{\rho}$ denotes the primitive ray generator of the ray $\rho$.
\end{thm}

Recall that a family of linear subspaces $\{V_\lambda \}_{\lambda \in \Lambda}$ of a finite dimensional vector space $V$ is said to form a distributive lattice if, there exists a basis $B$ of $V$ such that $B \cap V_\lambda$ is a basis of $V_\lambda$ for every $\lambda \in \Lambda$. If $X$ is nonsingular, then condition $({\bf C})$  of Theorem \ref{1} is equivalent to the following: for each $\sigma \in \Delta$, the collection of subspaces $\{E^{\rho}(i) \}_{\rho \in \sigma(1), i \in \Z}$ of $E$ forms a distributive lattice (see \cite[Remark $2.2.2$]{kly}). To see this, note that if the condition $({\bf C})$ is satisfied, for each $\sigma \in \Delta$, the eigen basis of the $\widehat{T}_{\sigma}$-grading $E=\oplus_{\chi \in \widehat{T}_{\sigma}}E^{[\sigma]}(\chi)$ serves as the required basis. Conversely, for $\sigma \in \Delta$, let $B_{\sigma}=\{e_1, \ldots, e_n \}$ be a basis of $E$ such that $E^{\rho}(i) \cap B_{\sigma}$ is a basis for $E^{\rho}(i)$. Then for each $e_j$ there exists an integer $n^j_{\rho}$ such that $e_j \in E^{\rho}(i) \cap B_{\sigma}$ for $i \leq n^j_{\rho}$ and $e_j \notin E^{\rho}(i) \cap B_{\sigma}$ for $i > n^j_{\rho}$. Since $\sigma$ is nonsingular cone, define  characters of $T_{\sigma}$ by $\chi_\textsubscript{\(j\)}:N_{\sigma} \rightarrow \Z$, $v_{\rho} \mapsto n^j_{\rho}$ for $\rho \in \sigma(1)$, which yields the required $\widehat{T}_{\sigma}$-grading $E=\oplus_{j=1}^n \text{Span}(e_j)$.

\begin{ex}[Filtrations for line bundles]\label{lb}{\rm
Let $\mathcal{L}=\mathcal{O}_X(D)$ be a toric line bundle on $X$ for some $T$-invariant Cartier divisor $D=\sum_{\rho \in \Delta(1)}a_{\rho}D_{\rho}$, $a_{\rho} \in \Z$. Then the associated filtrations $(L, \{L^{\rho}(i)\}_{\rho \in \Delta(1)})$ are given by:
\[ L^{\rho}(i) = \left\{ \begin{array}
{r@{\quad \quad}l}
L(=\C) & i \leqslant a_{\rho} \\ 

 0 & i > a_{\rho}
\end{array} \right. \] 

\noindent
Note that the filtrations define a function $f: \Delta(1) \rightarrow \Z$, given by $\rho \mapsto a_{\rho}$. Conversely when $X$ is nonsingular, any such function $f$ defines a toric line bundle whose associated filtrations are given by:
\[ L^{\rho}(i) = \left\{ \begin{array}
{r@{\quad \quad}l}
\C & i \leqslant f(\rho) \\ 

 0 & i >  f(\rho)
\end{array} \right. \]
}
\end{ex}

Note that for any two toric vector bundles $\mathcal{E}$ and $\mathcal{F}$ with the associated filtrations $(E, \{ E^{\rho}(i) \}_{\rho \in \Delta(1)})$ and $(F, \{ F^{\rho}(i) \}_{\rho \in \Delta(1)})$ respectively, the associated filtrations for $\mathcal{E \oplus F}$ and  $\mathcal{E \otimes F}$ are given by $(E \oplus F, \{E^{\rho}(i) \oplus F^{\rho}(i)\}_{\rho \in \Delta(1)})$ and $(E \otimes F, \{\sum_{s+ t=i} E^{\rho}(s) \otimes F^{\rho}(t) \}_{\rho \in \Delta(1)})$ respectively (see \cite[Lemma II.7.]{jgonza}). In particular, if $\mathcal{F}$ is a toric line bundle $\mathcal{O}_X(D)$, for some $T$-invariant Cartier divisor $D=\sum_{\rho \in \Delta(1)}a_{\rho}D_{\rho}$ with the associated  filtrations $(L, \{ L^{\rho}(i) \}_{\rho \in \Delta(1)})$ then the associated filtrations for $\mathcal{E} \otimes \mathcal{O}_X(D)$ is given by $(E , \{ E^{\rho}(i-a_{\rho}) \}_{\rho \in \Delta(1)} )$ (see \cite[Example II.9.]{jgonza}). 

\begin{ex}[Filtrations for tangent bundles]\label{tb}{\rm 
Let ${T}_X$ be the tangent bundle of the nonsingular toric variety $X=X(\Delta)$. Then the associated filtrations $(E, \{E^{\rho}(i)\}_{\rho \in \Delta(1)})$ are given by:\\
\[ E^{\rho}(i) = \left\{ \begin{array}
{r@{\quad \quad}l}
N_{\C} & i \leqslant 0 \\ 
\text{Span }( v_{\rho} ) & i=1 \\
 0 & i > 1
\end{array} \right. \]
}
\end{ex}

\noindent
We say that a toric vector bundle splits if it is equivariantly isomorphic to a direct sum of toric line bundles. Recall that, a toric vector bundle $\mathcal{E}$ over a nonsingular toric variety $X=X(\Delta)$ splits if and only if the filtrations $E^{\rho}(i)$, $\rho \in \Delta(1)$ generate a distributive lattice (see \cite[Corollary $2.2.3$]{kly}). Using this, we obtain a combinatorial criterion of splitting of the tangent bundle of a nonsingular toric variety in terms of the fan structure:

\begin{prop}\label{18}
Let $X=X(\Delta)$ be an \(n\)-dimensional nonsingular toric variety. Then ${T}_X$ splits if and only if for all $\sigma \in \Delta(n)$ and any ray $\rho \in \Delta(1) \setminus \sigma(1)$ the primitive ray generator $v_\rho$ of $\rho$ is the negative of some primitive ray generator of the cone $\sigma$.
\end{prop}

\noindent
{\bf Proof: }Suppose that ${T}_X$ splits. Then there exists a basis $B$ of $E~(=N_{\C})$ such that $B \cap E^\rho(1)$ is a basis of $E^\rho(1)=\text{Span}(v_\rho)$ for all $\rho \in \Delta(1)$. Let \(\sigma \in \Delta(n)\) and write \(\sigma =\text{Cone}(v_{\rho_1}, \ldots, v_{\rho_n})\). Then \(B\) must be of the form \(\{c_1 v_{\rho_1}, \ldots, c_n v_{\rho_n}\}\) for some non-zero scalars \(c_1, \ldots, c_n\). Let \(\rho \in \Delta(1) \setminus \sigma(1) \). Then since $B \cap E^\rho(1)$ is a basis of $E^\rho(1)$, \(v_{\rho}\) must be scalar multiple of one of \(c_1 v_{\rho_1}, \ldots, c_n v_{\rho_n}\). Since these are primitive vectors we must have \(v_{\rho}=-v_{\rho_j}\) for some \(j \in \{1, \ldots, n\}\), which concludes the forward direction.

Conversely, fix a cone $\sigma=\text{Cone }(v_{\rho_1}, \ldots, v_{\rho_n}) \in \Delta(n)$, then 
$\{v_{\rho_1}, \ldots, v_{\rho_n}\} \cap E^\rho(1)$ is a basis of $E^\rho(1)$ for all $\rho \in \Delta(1)$. Hence ${T}_X$ splits.
$\hfill\square$

\begin{cor}
Let $X=X(\Delta)$ be an \(n\)-dimensional nonsingular toric variety. If $|\Delta(1)| > 2n$, then ${T}_X$ does not split.
\end{cor}

\noindent
{\bf Proof: }On the contrary assume that \(T_X\) splits. Then fix \(\sigma \in \Delta(n) \) and let \(\rho \in \Delta(1) \setminus \sigma(1) \), then the primitive ray generator $v_\rho$ of $\rho$ is negative of some primitive ray generator of the cone $\sigma$ by Proposition \ref{18}. This implies $|\Delta(1)| \leq 2n$, which is a contradiction. Hence ${T}_X$ does not split.
$\hfill\square$

\subsection{Bott tower} 

\noindent
We now briefly recall the construction of Bott tower and some basic results related to it. For more details see \cite{grossbergkarshon} and \cite{civanyusuf}. Let $e_1, \ldots, e_r$ be the standard basis for $\mathbb{R}^r$, for any $r \in \N$. Given $\frac{n(n-1)}{2}$ integers $\{c_{i, j}\}_{\{1\leq i < j \leq n \}}$, the Bott tower of height $n$
$$M_n \rightarrow M_{n-1} \rightarrow \cdots \rightarrow M_2 \rightarrow M_1 \rightarrow M_0=\{ \text{point} \}$$
is constructed recursively as follows: Let $M_0$ be a point and $\xi_0$ be the trivial line bundle on $M_0$. Define $M_1$ to be $\mathbb{P}(\mathcal{O}_{M_0} \oplus \xi_0)=\mathbb{P}^1$. Let $\Delta_1 \subset \R$ be the fan of $M_1$ consisting of rays $v^1_1=e_1, v^1_2=-e_1$. Note that $\text{Pic}(M_1)=\Z [D_2]$, where $D_2$ is the invariant prime divisor corresponding to $v^1_2$. Let $\xi_1$ be the toric line bundle on $M_1$ whose associated filtrations are given by:
\begin{center}
$
 \xi_1^{v^1_1}(l) = \left\{ \begin{array}{ccc}

\C & l \leqslant 0 \\ 

 0 & l > 0
\end{array} \right. \ \  \ \ \text{and \ \ \ \ }
$
$
\xi_1^{v^1_2}(l) = \left\{ \begin{array}{ccc}

\C & l \leqslant c_{1, 2} \\ 

 0 & l > c_{1,  2} 
\end{array} \right. 
$
\end{center}

\noindent
Let $M_2:=\mathbb{P}(\mathcal{O}_{M_1} \oplus \xi_1)$ be the $\mathbb{P}^1$-bundle on $M_1$, where the fiber $\mathbb{P}^1$ has ray generators $e'_1, e'_0=-e'_1$. Consider the map $p : \R \rightarrow \R \oplus \R$, given by $v^1_1 \mapsto (v^1_1, 0)$ and $v^1_2 \mapsto (v^1_2, c_{1, 2}e'_1)$. Then the ray generators of $M_2$ are $v^2_1=e_1,v^2_2=e_2,v^2_3=-e_1 + c_{1, 2}e_2, v^2_4=-e_2$. The four maximal cones are $\text{Cone }(v^2_1, v^2_2), \text{Cone }(v^2_1, v^2_4), \text{Cone }(v^2_2, v^2_3)$ and $\text{Cone }(v^2_3, v^2_4)$. Repeating this process, at the $k$-th stage, we have a $k$-dimensional nonsingular projective toric variety $M_k$, whose fan structure is as follows: Let $v_i=e_i$ for $i=1, \ldots, k$, 
$v_{k+i}=-e_i+c_{i, i+1}e_{i+1}+\cdots +c_{i, k}e_k \text{ for }i=1, \ldots, k-1$
and $v_{2k}=-e_k$.\\
The fan $\Delta_k$ of $M_k$ is complete and consists of these $2k$ edges and $2^k$ maximal cones of dimension $k$ generated by these edges such that no cone contains both the edges $v_i$ and $v_{k+i}$ for $i=1, \ldots, k$.\\
Consider a line bundle $\xi_k$ on $M_k$ whose associated filtrations are given by:
\begin{center}
$
\xi_k^{v_i}(l) = \left\{ \begin{array}{ccc}

\C & l \leqslant 0 \\ 

 0 & l > 0
\end{array} \right. \ \ 
$ 
and \ \ \ \
$ 
\xi_k^{v_{k+i}}(l) = \left\{ \begin{array}{ccc}

\C & l \leqslant c_{i, k+1} \\ 

 0 & l > c_{i, k+1} 
\end{array} \right. 
$ 
for $i=1, \ldots, k$.
\end{center}

\noindent
Then $M_{k+1}:=\mathbb{P}(\mathcal{O}_{M_k} \oplus \xi_k)$. The integers $c_{i, j}$'s are called Bott numbers and they are arranged in a $n \times n$ upper triangular matrix with $1$'s in the diagonal, called Bott matrix.

\noindent
Let $D_i:=D_{v_i}$ denote the invariant prime divisor corresponding to an edge $v_i$ for \(i=1, \ldots, 2k\). Note that for any \(D \in \text{Pic}(M_k)\), \(D=\sum\limits_{i=1}^{2k} a_i D_i \), where \(a_1, \ldots, a_{2k}\) are integers (see \cite[Theorem \(4.2.1\)]{Cox}). We also have the following relations:
\begin{center}
$D_i \sim_{\text{lin}} D_{k+i}-c_{1, i}D_{k+1}-\cdots-c_{i-1, i}D_{k+i-1}$ for $i=1, \ldots, k$.
\end{center}
Hence any divisor $D$ on $M_k$ can be written as \(D=\sum\limits_{i=1}^{k}b_iD_{k+i},\) where $b_i \in \Z$ for $i=1, \ldots, k$. Now we show that \(D_{k+1}, \ldots, D_{2k}\) are linearly independent. Suppose 
\begin{equation}\label{li}
c_1 D_{k+1} + \cdots + c_{k} D_{2k}=0 
\end{equation}
for some scalars \(c_1, \ldots, c_k \in \Z \). Consider the cones \(\tau_i=\text{Cone}(v_1, \ldots, \widehat{v}_i, \ldots, v_k)\) for $i=1, \ldots, k$. Since $D_{k+l} \cdot V(\tau_i)=0$ for $l \neq i$ and $D_{k+i} \cdot V(\tau_i)=1$ (see \cite[Corollary $6.4.3$]{Cox}), taking intersection product of \eqref{li} with \(V(\tau_i)\) we see that  \(c_i=0\) for $i=1, \ldots, k$. Hence the Picard group of the Bott tower is given by $\text{Pic}(M_k)=\Z [D_{k+1}] \oplus \ldots \oplus \Z [D_{2k}]$.

\vspace{1cm}
\noindent 
The key observation regarding the Bott numbers is the following:

\begin{thm}\label{12} 
The integers $\{c_{i, j}\}_{\{1\leq i < j \leq n \}}$ can be assumed to be non-negative.
\end{thm}

\noindent
{\bf Proof:} We prove it by induction on $n$.\\
For $n=2$ the corresponding integer is $c_{1, 2}$ and $M_2=\mathcal{H}_{c_{1, 2}}$, Hirzebruch surface and it is well known that $c_{1, 2}$ can be assumed to be non-negative. By induction hypothesis assume that $c_{i, j}$'s are non-negative for $1\leq i < j \leq k $.\\
Consider a line bundle $\eta_k$ on $M_k$ with associated filtrations:
\begin{center}
$
 \eta_k^{v_i}(l) = \left\{ \begin{array}{ccc}

\C & l \leqslant 0 \\ 

 0 & l > 0
\end{array} \right. \ \ \ 
$ 
and \ \ \ 
$
 \eta_k^{v_{k+i}}(l) = \left\{ \begin{array}{ccc}

\C & l \leqslant |c_{i, k+1}| \\ 

 0 & l > |c_{i, k+1}| 
\end{array} \right. 
$  
for $i=1, \ldots, k$.
\end{center}

\noindent
We show that there is a line bundle $L_k$ on $M_k$ such that $(\mathcal{O}_{M_k} \oplus \xi_k)\otimes L_k \cong \mathcal{O}_{M_k} \oplus \eta_k$, where $\xi_k$ is as before. Let $S=\{v_{k+i} | c_{i, k+1} < 0 \}$ and the filtrations associated to $L_k$ be as follows:\\
If $\alpha \notin S$, $ L_k^{\alpha}(l) = \left\{ \begin{array}{ccc}

\C & l \leqslant 0 \\ 

 0 & l > 0
\end{array} \right. \ \
$  
and if $\alpha \in S$, say $\alpha=v_{k+i}$ ,
$L_k^{\alpha}(l) = \left\{ \begin{array}{ccc}

\C & l \leqslant -c_{i, k+1}  \\ 

 0 & l > -c_{i, k+1}
\end{array} \right. 
$ \\

\noindent
Note that the associated filtrations for $\mathcal{O}_{M_k} \oplus \eta_k$ are given by:
\begin{center}
$
(\mathcal{O}_{M_k} \oplus \eta_k)^{v_i}(l) = \left\{ \begin{array}{ccc}

\C \oplus \C & l \leqslant 0 \\ 

 0 & l > 0
\end{array} \right. 
$
\end{center}
\begin{center}
when $\alpha=v_{k+i} \notin S $,
$(\mathcal{O}_{M_k} \oplus \eta_k)^{\alpha}(l) = \left\{ \begin{array}{ccc}

\C \oplus \C & l \leqslant 0 \\ 
\C    &   0 < l \leq c_{i, k+1}            \\
 0 & l > c_{i, k+1}
\end{array} \right. 
$
\end{center}
\begin{center}
when $\alpha \in S$, say $\alpha=v_{k+i}$, 
$ (\mathcal{O}_{M_k} \oplus \eta_k)^{\alpha}(l) = \left\{ \begin{array}{ccc}

\C \oplus \C & l \leqslant 0 \\ 
\C    &   0 < l \leq -c_{i, k+1}           \\
 0 & l > -c_{i, k+1}
\end{array} \right.
$

\end{center} 

Similarly, the filtrations associated to the vector bundle $( \mathcal{O}_{M_k} \oplus \xi_k)$ are given by:
\begin{center}
	$
	(\mathcal{O}_{M_k} \oplus \xi_k)^{v_i}(l) = \left\{ \begin{array}{ccc}
	
	\C \oplus \C & l \leqslant 0 \\ 
	
	0 & l > 0
\end{array} \right. 
$
\end{center}
\begin{center}
when $\alpha=v_{k+i} \notin S $,
$(\mathcal{O}_{M_k} \oplus \xi_k)^{\alpha}(l) = \left\{ \begin{array}{ccc}
	
	\C \oplus \C & l \leqslant 0 \\ 
	\C    &   0 < l \leq c_{i, k+1}            \\
	0 & l > c_{i, k+1}
\end{array} \right. 
$
\end{center}
\begin{center}
when $\alpha \in S$, say $\alpha=v_{k+i}$, 
$ (\mathcal{O}_{M_k} \oplus \xi_k)^{\alpha}(l) = \left\{ \begin{array}{ccc}
	
	\C \oplus \C & l \leqslant c_{i, k+1} \\ 
	\C    &   c_{i, k+1} < l \leq 0           \\
	0 & l > 0
\end{array} \right.
$

\end{center} 

\noindent
By similar computations, the filtrations associated to the vector bundle $( \mathcal{O}_{M_k} \oplus \xi_k)\otimes L_k$ are same as the filtrations associated to $\mathcal{O}_{M_k} \oplus \eta_k$. Hence by Theorem \ref{1} these two are $T$-equivariantly isomorphic. Thus $M_{k+1}:=\mathbb{P}(\mathcal{O}_{M_k} \oplus \xi_k)$ and $\mathbb{P}(\mathcal{O}_{M_k} \oplus \eta_k)$ are isomorphic as abstract varieties, hence they are isomorphic as toric varieties also (see \cite[Proposition \(3.1\)]{masuda}, \cite[Theorem $4.1$]{berchtold2003lifting}).
$\hfill\square$\\

\noindent
Let $T_{M_k}$ be the tangent bundle of $M_k$ with associated filtrations $(E, \{ E^\rho(i)\}_{\rho \in \Delta_k, i \in \Z})$ (see Example \ref{tb}).

\begin{prop}\label{15}
$T_{M_k}$ splits if and only if $c_{i, j}=0$ for all $ i \text{ and } j$.
\end{prop}

\noindent
{\bf Proof:} Suppose that $T_{M_k}$ splits. Consider the maximal cone \( \sigma= \text{Cone}(v_1, \ldots, v_k) \). Then by Proposition \ref{18}, $v_{k+i}$ has to be negative of some element of \(\{v_1, \ldots, v_k\} \) for \(i=1, \ldots, k \). Therefore $c_{i, j}=0$ for all \(i\) and \(j\).

Conversely, let $c_{i, j}=0$ for all $i$ and $j$. Then the choice of the basis $B= \{ e_1, \ldots, e_k\}$ of $\C^k$ works.$\hfill\square$ 

\begin{rmk}{\rm
In general, consider the toric line bundles $\mathcal{L}_0, \mathcal{L}_1, \ldots, \mathcal{L}_m$ over a nonsingular toric variety $X=X(\Delta)$. Then $X'= \mathbb{P}(\mathcal{L}_0 \oplus \mathcal{L}_1 \oplus \ldots \oplus \mathcal{L}_m)$ is also a nonsingular toric variety with fan $\Delta'$ given as follows: Let $\mathcal{L}_0, \mathcal{L}_1, \ldots, \mathcal{L}_m$ correspond to $\Delta$-linear support functions $h_0, h_1, \ldots, h_m$ respectively. Choose a $\Z$-basis $\{e_1, \ldots, e_m\}$ of $\R^m$ and let $e_0=-e_1- \ldots - e_m$. Consider the $\R$-linear map $\Phi: \R^n \rightarrow \R^n \oplus \R^m$ which sends $y \mapsto (y, -\sum_{j=0}^m h_j(y)e_j)$. Now let $\tilde{\sigma}_i = \text{Cone}(e_0, \ldots, \widehat{e_i}, \ldots, e_m)$ for each $0\leq i \leq m$ (here by $\widehat{e}_i$ we mean that $e_i$ is omitted from the relevant collection). Let $\tilde{\Delta}$ be the fan in $\R^m$ generated by $\tilde{\sigma}_i$ for $0\leq i \leq m$. Then $\Delta'=\{\Phi(\sigma) + \tilde{\sigma} : \sigma \in \Delta, \tilde{\sigma}\in \tilde{\Delta}\}$.  

\begin{enumerate}

\item For $m \geq 2$, $T_{X'}$ does not split. To see this, consider the maximal cone $\sigma'= \Phi(\sigma) + \text{Cone}(e_0, \ldots, \widehat{e_i}, \ldots, e_m)$, where $\sigma \in \Delta$ is a maximal cone. The ray generated by $(0, e_i)$ is not in $\sigma'$ and  $(0, e_i)$ is not the negative of any of the primitive ray generators of $\sigma'$.

\item For $m=1$, if $T_{X'}$ splits then $T_X$ also splits.

\end{enumerate}
}
\end{rmk}

\subsection{s-jet ampleness}

\noindent
Let us recall the notion of $s$-jet ampleness, see \cite{rider} for details. Let $X$ be a nonsingular algebraic variety. For a point $x \in X$, let $\mathfrak{m}_x$ be the maximal ideal sheaf of $x$ in $X$. 

\begin{defn}
A line bundle $\mathcal{L}$ on $X$ is said to be $s$-jet ample, if for every finite collection of points $x_1, \ldots, x_r$ the restriction
map $$H^0(X; \mathcal{L}) \longrightarrow H^0\left(X; \mathcal{L} \otimes \frac{\mathcal{O}_X}{\mathfrak{m}_{x_1}^{k_1} \otimes \cdots \otimes \mathfrak{m}_{x_r}^{k_r}} \right)$$
is surjective, where $k_1, \ldots, k_r$ are positive integers and $\sum_{i=1}^r k_i = s + 1$. 
\end{defn}

Note that $\mathcal{L}$ is $0$-jet ample if and only if $\mathcal{L}$ is globally generated and $\mathcal{L}$ is $1$-jet ample if and only if $\mathcal{L}$ is very ample (see \cite[Section $2.3$]{jetample}). On nonsingular complete toric varieties the notion of ample and very ample line bundles coincide (see \cite[Theorem $6.1.15$]{Cox}). Let us recall the following generalization of Toric-Nakai criterion for $s$-jet ampleness.

\begin{thm} \cite[Proposition $3.5$, Theorem $4.2$]{DiRocco}\label{11} Let $\mathcal{L}$ be a line bundle on a nonsingular toric variety $X=X(\Delta)$. Then the following are equivalent:
\begin{itemize}
\item[(1)] $\mathcal{L}$ is $s$-jet ample.
\item[(2)] $\mathcal{L} \cdot C \geq s$, for any $T$-invariant curve C.
\end{itemize}
\end{thm}

\section{Positivity of line bundles on Bott towers}

\subsection{s-jet ampleness on Bott tower}

\noindent
By the virtue of Theorem \ref{12}, we can assume that the Bott numbers are non-negative without loss of generality. We use Theorem \ref{11} to give a necessary and sufficient criterion for $s$-jet ampleness of line bundles on Bott tower.

\begin{thm}\label{16}
Let $D=\sum_{i=1}^{k}a_iD_{k+i}$ be a $T$-invariant Cartier divisor on $M_k$. Then $D$ is  $s$-jet ample if and only if $a_i  \geq \ s$ for all $i=1, \ldots, k$.
\end{thm}
\noindent
{\bf Proof:} Let $D$ be $s$-jet ample. For $i=1, \ldots, k$, consider the cones $\tau_i=\text{Cone}(v_1, \ldots, \widehat{v}_i, \ldots, v_k)$. Now $D \cdot V(\tau_i)  \geq s$ for $i=1, \ldots, k$ by Theorem \ref{11}. Note that $D \cdot V(\tau_i)=a_i$ since $D_{k+l} \cdot V(\tau_i)=0$ for $l \neq i$ and $D_{k+i} \cdot V(\tau_i)=1$ (\cite[Corollary $6.4.3$]{Cox}).

Conversely, let $a_i  \geq \ s$ for all $i=1, \ldots, k$. Again  by Theorem \ref{11}, it suffices to show that $D \cdot V(\tau) \geq s$ for any $k-1$ dimensional cone $\tau \in \Delta_k$. Note that any such cone is of the form 
\begin{center}
$\tau=\text{Cone}(v_1, \ldots, \widehat{v}_{i_\textsubscript{1}}, \ldots, \widehat{v}_{i_\textsubscript{\(r\)}}, \ldots, v_k,v_{k+i_\textsubscript{1}}, v_{k+i_\textsubscript{\(2\)}}, \ldots, \widehat{v}_{k+i_\textsubscript{\(j\)}}, \ldots, v_{k+i_\textsubscript{\(r\)}})$ 
\end{center}
for some $j=1, \ldots, r$, where \(r\) varies from \(1\) to \(k\). We have\\
$
 D_{k+m} \cdot V(\tau) = \left\{ \begin{array}{ccc}

1 & m=i_j   \\
0 & m \neq i_1, \ldots, i_r \\
\end{array} \right.  \ \
$
and \ \ \ 
$
  D_{k+m} \cdot V(\tau) = \left\{ \begin{array}{ccc}

 0 &   m=i_l, 1 \leq l <j ,\\
 b_{i_\textsubscript{\(l\)}} & m =i_l, j < l \leq r ,\\
\end{array} \right. \\ 
$

\noindent
where $b_{i_l}$'s are such that the following wall relation is satisfied: 
\begin{center}
$v_{i_\textsubscript{\(j\)}}+d_1v_1+\cdots+\widehat{v}_{i_\textsubscript{\(1\)}}+\cdots+\widehat{v}_{i_\textsubscript{\(r\)}}+\cdots+d_kv_k+v_{k+i_\textsubscript{\(j\)}}+\sum_{t=j+1}^r b_{i_\textsubscript{\(t\)}}v_{k+i_\textsubscript{\(t\)}}=0$
\end{center}
where $d_1, \ldots, d_k$ are integers(\cite[Corollary $6.4.3$ and Proposition $6.4.4$]{Cox}). Then we have, 
\begin{center}
$b_{i_\textsubscript{\(j+1\)}}=c_{i_\textsubscript{\(j\)}, i_\textsubscript{\(j+1\)}}$; $b_{i_\textsubscript{\(j+2\)}}=c_{i_\textsubscript{\(j+1\)}, i_\textsubscript{\(j+2\)}} b_{i_\textsubscript{\(j+1\)}}+c_{i_\textsubscript{\(j\)}, i_\textsubscript{\(j+2\)}}$;\\
  $\vdots$\\
  $b_{i_\textsubscript{\(r\)}}=c_{i_\textsubscript{\(j+1\)}, i_\textsubscript{\(r\)}} b_{i_\textsubscript{\(j+1\)}} + c_{i_\textsubscript{\(j+2\)}, i_\textsubscript{\(r\)}} b_{i_\textsubscript{\(j+2\)}}+ \cdots+c_{i_\textsubscript{\(r-1\)}, i_\textsubscript{\(r\)}} b_{i_\textsubscript{\(r-1\)}} + c_{i_\textsubscript{\(j\)},\ i_\textsubscript{\(r\)}} $.
\end{center}
Hence $D \cdot V(\tau)=a_{i_j}+\sum_{l=j+1}^r a_{i_l}  b_{i_l} \geq s$ since $c_{i, j} \geq 0$ by Theorem \ref{12}.
$\hfill\square$\\

\noindent
We have the following applications of the above theorem.

\begin{cor}\label{2}
The Cartier divisor $D=\sum_{i=1}^{k}a_iD_{k+i}$ on $M_k$ is ample (respectively, nef) if and only if $a_i >  (\text{respectively, } \geq) \ 0$ for all $i=1, \ldots, k$. The anticanonical divisor $-K_{M_k}= D_1 + \ldots + D_{2k}$ is at most $2$-jet ample. Moreover, it is $2$-jet ample if and only if $c_{i, j}=0$ for all $i$ and $j$. Further, $M_k$ is Fano if and only if $\sum_{j=i+1}^k c_{i, j} \leq 1, i=1, \ldots, k-1$. Finally, the cotangent bundle on $M_k$ is never ample since $ K_{ M_k}$ is not ample.
\end{cor}

\begin{rmk}{\rm
Let $D=\sum_{i=1}^{k}a_iD_{k+i}$ be an ample $T$-invariant Cartier divisor on $M_k$. Then by the above ampleness criterion, $D-D_j$ is ample if and only if 
\begin{itemize}
\item $a_j >1$ if $1 \leq j \leq k$ and
\item $a_{j-k}>1$ if $k+1 \leq j \leq 2k$.
\end{itemize}
For $j <j'$, $D-D_j-D_{j'}$ is not nef if and only if $j'=k+j$ and $a_j=1$. In fact, $D-D_j-D_{j'}$ is ample if and only if 
\begin{itemize}
\item $a_j > 2$, when $j'=k+j$.
\item $b_j > 1$ and $b_{j'} > 1$, otherwise where 
\[ b_l = \left\{ \begin{array}
{r@{\quad \quad}l}
a_l & l \leqslant k \\ 

 a_{l-k} & l > k
\end{array} \right. \]   

\end{itemize}
This generalizes Mustata's criterion (see \cite{van}) in case of Bott tower.
}
\end{rmk}

\noindent
{\bf Generalization of ampleness criterion for toric vector bundles on Bott tower:}

Let us recall the notion of semi-stability of a nonzero torsion-free coherent sheaf $\mathcal{E}$ on $X$, where $X$ is a normal projective variety of dimension $n$ with an ample divisor $H$ on $X$. The slope of $\mathcal{E}$ is defined to be $\mu(\mathcal{E}):=\frac{\text{deg }(\mathcal{E})}{\text{rank }\mathcal{E}}$, where $\text{deg }(\mathcal{E})=c_1(\mathcal{E}) \cdot H^{n-1}$. The sheaf $\mathcal{E}$ is then said to be semi-stable if $\mu(\mathcal{F}) \leq \mu(\mathcal{E})$ for all proper nonzero subsheaves $\mathcal{F}\subsetneq \mathcal{E}$.

 If $\mathcal{E}$ is a semi-stable vector bundle on a smooth projective curve over an algebraically closed field of characteristic zero, it is enough to consider its determinant bundle to check for ampleness or nefness of the vector bundle (see  \cite[Theorem $3.2.7$]{huybrechts2010geometry}). The following proposition generalizes this result for a class of semi-stable toric vector bundles over any nonsingular projective complex toric variety.

\begin{prop}\label{7} 
Let $\mathcal{E}$ be a toric vector bundle on a nonsingular projective toric variety $X=X(\Delta)$ of rank \(m\). Assume that $\mathcal{E}$ is semi-stable and discriminant $\Delta(\mathcal{E})=0$, where the discriminant of $\mathcal{E}$ is defined to be the characteristic class $\Delta(\mathcal{E}):=c_2(\mathcal{E}) - \frac{m-1}{2m} c_1(\mathcal{E})^2$. Then $\mathcal{E}$ is ample (respectively, nef) if and only if $\text{det}(\mathcal{E})$ is ample $($respectively, nef$)$.
\end{prop}

\noindent
{\bf Proof:} By \cite[Theorem $2.1$]{hering2010positivity}, $\mathcal{E}$ is ample (respectively, nef) if and only if $\mathcal{E}|_{V(\tau)}$ is ample (respectively, nef) for all walls $\tau$ (here wall means codimension $1$ cone). Since $\mathcal{E}$ is semi-stable discriminant zero, we have \(\mathcal{E}|_{V(\tau)}\) is semi-stable for all wall $\tau$ by \cite[Theorem $2.5$]{bruzzo2013restricting}. Hence $\mathcal{E}|_{V(\tau)}$ is ample (respectively, nef) for all walls $\tau$ if and only if $\text{deg}(\mathcal{E}|_{V(\tau)})> (\text{respectively, } \geq) 0$ for all walls $\tau$ by \cite[Theorem $3.2.7$]{huybrechts2010geometry}. Now $\text{deg}(\mathcal{E}|_{V(\tau)})=\text{det}(\mathcal{E}) \cdot V(\tau)$. Hence we have  $\mathcal{E}$ is ample (respectively, nef) if and only if $\text{det}(\mathcal{E}) \cdot V(\tau) > (\text{respectively, }\geq) 0$  for all walls $\tau$ if and only if $\text{det}(\mathcal{E})$ is ample (respectively, nef) by Toric Nakai criterion \cite[Theorem $2.18$]{Oda}. $\hfill\square$

\subsection{Criterion for nef and big line bundle on Bott tower}
In this subsection we give a criterion for nef and big line bundle on Bott tower. A characterization for big divisors on an irreducible projective variety is given in \cite[Corollary $2.2.7$]{lazarsfeld2004positivity}. The following proposition is an equivariant version of the above characterization for big divisors on nonsingular projective toric varieties.

\begin{prop}[Characterization of $T$-invariant big divisors]\label{3}
	Let $X=X(\Delta)$ be a nonsingular projective toric variety. Let $D$ be a $T$-invariant Cartier divisor on X, then the following are equivalent:
	\begin{itemize}
		\item[(i)]$D$ is big.
		\item[(ii)]For any $T$-invariant ample integer divisor $A$ on $X$, there exists a positive integer $n$ and $T$-invariant effective divisor $N$ on $X$, such that $nD \sim_{\text{lin}} A+N$.
		\item[(iii)] same as in $(ii)$ for some $T$-invariant ample divisor $A$ on $X$.
		\item[(iv)]There exists $T$-invariant ample divisor $A$ on $X$, a positive integer $n$  and $T$-invariant effective divisor $N$ on $X$, such that $nD \equiv_{\text{num}} A+N$.  
	\end{itemize}
\end{prop}

\noindent
{\bf Proof:} Assume $D$ is big. Since $A$ is $T$-invariant ample, there exists a positive integer $r$ such that both $rA$ and $(r+1)A$ are globally generated. So we can choose nonzero global sections $\chi^{m_1} \in H^0(X, \mathcal{O}_X(rA))$ and  $\chi^{m_2} \in H^0(X, \mathcal{O}_X((r+1)A))$. Then $rA \sim_{\text{lin}}(\chi^{m_1})_0$ and $(r+1)A \sim_{\text{lin}}(\chi^{m_2})_0$, where $H_r:=(\chi^{m_1})_0, \  H_{r+1}:=(\chi^{m_2})_0$ are effective divisors (\cite[Proposition II $7.7$]{hartshorne2013algebraic}). In fact, both the divisors $H_r$ and $H_{r+1}$ are also $T$-invariant. To see this, let \(\{(U_{\sigma}, \chi^{-m_{\sigma}})\}_{\sigma \in \Delta}\) be a local data for \(rA\), i.e. \(rA|_{U_{\sigma}}=\text{div}(\chi^{-m_{\sigma}})|_{U_{\sigma}}\) (see \cite[Theorem 4.2.8]{Cox} ). Hence we have the isomorphism $\phi_\sigma: \mathcal{O}_X(rA)\mid_{U_\sigma} \simeq \mathcal{O}_{U_\sigma}$ given by multiplication by $\chi^{-m_\sigma}$. Then the Cartier data for $H_r$ is given by $\{m_1-m_\sigma\}_{\sigma \in \Delta}$, which shows that $H_r$ is $T$-invariant and similarly $H_{r+1}$ is also $T$-invariant.
Now by Kodaira's Lemma (\cite[Proposition $2.2.6$]{lazarsfeld14500positivity}) we have $H^0(X, \mathcal{O}_X(nD-H_{r+1})) \neq 0$ for some large positive integer $n$. Since $nD-H_{r+1}$ is $T$-invariant, arguing as above can choose nonzero global section $\chi^{m_3} \in H^0(X, \mathcal{O}_X(nD-H_{r+1}))$. Then $nD-H_{r+1} \sim_{\text{lin}} nD-H_{r+1}+\text{div}(\chi^{m_3} ) \geq 0$. Take $N'= nD-H_{r+1}+\text{div}(\chi^{m_3} )$ which is effective and $T$-invariant. Now  $nD-H_{r+1} \sim_{\text{lin}} N'$, which implies $nD \sim_{\text{lin}} H_{r+1}+N' \sim_{\text{lin}} (r+1)A +N' \sim_{\text{lin}} A +H_r+N'$. Taking $N=H_r+N'$, which is $T$-invariant effective we get  $nD \sim_{\text{lin}} A+N$. \\
\noindent
$(ii)\Rightarrow (iii)\Rightarrow (iv)$ Obvious.\\
\noindent
$(iv)\Rightarrow (i)$ Follows from the same argument as in (\cite[Corollary $2.2.7$]{lazarsfeld2004positivity}).
$\hfill\square$

\begin{cor}\label{4}
Let $X=X(\Delta)$ be a nonsingular projective toric variety. Let $D$ be a $T$-invariant Cartier divisor on $X$, then $D$ is nef and big if and only if there exists $T$-invariant effective divisor $N$ on $X$, such that $D-\frac{1}{n} N$ is ample for $n\gg 0$.
\end{cor}

\noindent
{\bf Proof:} Let $D$ be nef and big. Then by Proposition \ref{3} there exists a positive integer $k$, $T$-invariant ample Cartier divisor $A$ and $T$-invariant effective Cartier divisor $N$ on $X$, such that $kD \equiv_{\text{num}} A+N$. Take $n > k$, then $nD\equiv_{\text{num}} (n-k)D+A+N$. Since $D$ is nef and $A$ is ample, by \cite[Corollary $1.4.10$, Theorem $1.2.23$, Theorem $1.4.9$]{lazarsfeld14500positivity} $(n-k)D+A$ is ample. So $nD-N $ is ample, which implies $D-\frac{1}{n} N$ is ample for $n\gg 0$.

Converse direction follows by the same argument as in \cite[Example 2.2.19]{lazarsfeld14500positivity}.
$\hfill\square$

\begin{rmk}\label{19}\rm
Using Corollary \ref{4}, for a $T$-invariant divisor $D=\sum_{i=1}^{k}a_iD_{k+i}$ on $M_k$ to be big one must have $a_k>0$. Hence the canonical divisor $K_{M_k}$ is never big.

\end{rmk}

\noindent
We now give a necessary and sufficient criterion for nef and big divisor on Bott tower.
 
\begin{thm}\label{5}
	Let $D=\sum_{i=1}^{k}a_iD_{k+i}$ be a $T$-invariant Cartier divisor on $M_k$. Then $D$ is nef and big if and only if \(a_i'\)s satisfies the following conditions:
	\begin{itemize}
		\item $a_k>0$, $a_i \geq 0$ for $1 \leq i < k$ and
		\item if \(c_{i, j}=0\) for all \(j>i\) for some \(i\), then the corresponding $a_i$ must be positive.
	\end{itemize}
\end{thm}

\noindent
{\bf Proof:} Let $D$ be nef and big. Then by Corollary \ref{2} and Remark \ref{19} we have $a_i \geq 0$ for $1 \leq i \leq k-1$ and $a_k>0$. Again Corollary \ref{4} implies that, there exists $T$-invariant effective divisor $N=\sum_{i=1}^{2k}b_iD_{k}$ on $M_k$ $(b_i \geq 0 \text{ for all } i)$, such that $D-\frac{1}{n} N$ is ample for $n\gg 0$. Write $N \sim_{\text{lin}} \sum_{i=1}^{k}c_iD_{k+i} $, where $c_k=b_k + b_{2k}$ and $c_i=b_i+b_{k+i}- \sum_{j=i+1}^kc_{i, j}b_j$, for $i=1, \ldots, k-1$. Take $d_i=a_i-\frac{1}{n} c_i$ and write $D-\frac{1}{n} N=\sum_{i=1}^{k}d_iD_{k+i}$. Now if \(c_{i, j}=0\) for all \(j>i\) for some \(i\), we have \(d_i=a_i-\frac{1}{n} (b_i+b_{k+i}) \). By Corollary \ref{2}, $d_i$ must be positive, which is possible only if $a_i$ is positive. 

Conversely, by Corollary \ref{2} we have $D$ is nef. Consider the $T$-invariant effective divisor $N=\sum_{i=1}^{2k}b_iD_i$ on $M_k$, where $\{ b_1 , b_2 , \ldots , b_k\}$ is a strictly increasing sequence of positive integers and $b_{k+i}=0$ for $i=1, \ldots, k$. Write $N \sim_{\text{lin}} \sum_{i=1}^{k}c_iD_{k+i} $, with $c_k=b_k $ and $c_i=b_i- \sum_{j=i+1}^kc_{i, j}b_j$, for $i=1, \ldots, k-1$. By Corollary \ref{4} it suffices to show that $D-\frac{1}{n} N=\sum_{i=1}^k(a_i- \frac{1}{n}c_i) D_i$ is ample for $n\gg 0$. Since $a_k \geq 1$, for some positive integer $n_0$ we have $a_k- \frac{1}{n}c_k>0$ whenever $n > n_0$. Now if \(i\) is such that \(c_{i, j}=0\) for all \(j>i\), then  $ a_i - \frac{1}{n}c_i=a_i-\frac{1}{n} b_i$ and this is positive whenever $n > n_i$ for some large positive integer $n_i$ since in this case $a_i$ is positive. Now if \(i\) is such that \(c_{i, j'} \neq 0\) for some $j'>i$. Then $a_i - \frac{1}{n}c_i=a_i- \frac{1}{n}(b_i- \sum_{j=i+1}^k c_{i, j}b_j )=a_i- \frac{1}{n}(b_i-c_{i, j'}b_{j'}) + \frac{1}{n} \left( \sum_{j=i+1,j \neq j'}^k c_{i, j}b_j \right)> 0$ for any \(n\) since $b_i-c_{i, j'}b_{j'} < 0$ and $\sum_{j=i+1,j \neq j'}^k c_{i, j}b_j \geq 0$ by the choice of $b_i$'s. Hence $a_i - \frac{1}{n}c_i > 0$ for all \(i=1, \ldots, k\) and for $n\gg 0$. So $D-\frac{1}{n} N$ is ample for $n\gg 0$. This completes the proof.
$\hfill\square$ 	 

Recall that a nonsingular projective algebraic variety is said to be weak Fano if its anticanonical divisor is nef and big. As an immediate corollary we obtain the following.

\begin{cor}\label{17}
$M_k$ is weak Fano if and only if $\sum_{j=i+1}^k c_{i, j} \leq 2, i=1, \ldots, k-1$.
\end{cor}

\begin{cor}
Every nef and big divisor on $M_k$ is ample if and only if $c_{i, j}=0$ for all $i$ and $j$. 
\end{cor}
\noindent
{\bf Proof:} Let $c_{i, j} \neq 0$ for some $i, \ j$. Then consider $D= D_{k+1}+ \ldots + D_{k+i-1}+D_{k+i+1}+\ldots +D_{2k}$. Hence $D$ is nef and big by Theorem \ref{5}, but not ample by Corollary \ref{2}.

Converse is immediate from Theorem \ref{5} and Corollary \ref{2}.
$\hfill\square$

\section{Toric subbundles of a toric vector bundle}

Let \(\mathcal{E}\) be a toric vector bundle on a nonsingular toric variety \(X\). Recall that a toric vector bundle \(\mathcal{F}\) on \(X\) is said to be a toric subbundle of \(\mathcal{E}\) if it is a sub vector bundle in the usual sense $($i.e., \( \mathcal{F} \subseteq \mathcal{E} \)  and for each point \(x \in X\), the fiber \(\mathcal{F}(x)\) is a vector subspace of \(\mathcal{E}(x)\)$)$ and the above inclusion \( \mathcal{F} \subseteq \mathcal{E} \) of the corresponding total spaces is \(T\)-equivariant such that the following diagram commutes for all \(x \in X\) and for all \(t \in T\).

\begin{center}
	\begin{tikzpicture}[description/.style={fill=white,inner sep=2pt}]
	\matrix (m) [matrix of math nodes, row sep=3em,
	column sep=1.5em, text height=1.5ex, text depth=0.25ex]
	{  \mathcal{F}(x)  & & \mathcal{E}(x) \\
		\mathcal{F}(t \cdot x) & &  \mathcal{E}(t \cdot x)  \\ };
	\path[right hook->]  (m-1-1) edge node[auto] {}(m-1-3);
	\path[right hook->] (m-2-1) edge node[below] {}(m-2-3);
	\path[->] (m-1-1) edge node[auto] {$t \cdot $}(m-2-1);
	\path[->] (m-1-3) edge node[auto] {$ t \cdot $} (m-2-3);
	
	\end{tikzpicture}
\end{center}
Note that if \(\mathcal{F}\) on \(X\) is a toric subbundle of \(\mathcal{E}\), then \(\frac{\mathcal{E}}{\mathcal{F}} \) is a vector bundle (See \cite[Section 1.7]{lepotier}).

\subsection{Characterization of toric subbundles of a toric vector bundle}\label{toricsub}

\noindent
Here we follow the notations introduced in Section 2. Let $X=X(\Delta)$ be a nonsingular toric variety and $\mathcal{E}$ be a toric vector bundle on $X$ with associated filtrations $\left( E, \{E^{\rho}(i) \}_{\rho \in \Delta(1)} \right)$. Let $\mathcal{F}$ be a toric subbundle with associated filtrations $\left( F, \{F^{\rho}(i) \}_{\rho \in \Delta(1)} \right)$. We characterize the filtration data for $\mathcal{F}$ as follows:
\begin{prop}\label{6}
	Let $ \mathcal{E}$ be a toric vector bundle on a nonsingular toric variety $X$ with associated filtrations $\left( E, \{E^{\rho}(i) \}_{\rho \in \Delta(1)} \right)$. Then the toric subbundles of $\mathcal{E}$ are in one-to-one correspondence with subfiltrations $\left( F, \{F^{\rho}(i) \}_{\rho \in \Delta(1)} \right)$ of $\left( E, \{E^{\rho}(i) \}_{\rho \in \Delta(1)} \right)$, where $F \subseteq E$, $F^{\rho}(i)=F \cap E^{\rho}(i)$ and the collection of subspaces $\{F, \{E^{\rho}(i) \}_{\rho \in \sigma(1)} \}$ of $E$ forms a distributive lattice for all $\sigma \in \Delta$.
\end{prop}
\noindent
{\bf Proof:} Let $ \mathcal{E}$ be a toric vector bundle on $X$ and $\mathcal{F}$ be a toric subbundle of $\mathcal{E}$. Fix $\sigma \in \Delta$ and consider the evaluation map \(ev_{x_{\sigma}} : \Gamma(U_{\sigma}, \mathcal{E}) \rightarrow \mathcal{E}(x_{\sigma})\), given by evaluating the sections at the distinguished point \(x_{\sigma}\). Since \(\mathcal{F}\) is a toric subbundle of \(\mathcal{E}\) we have the following commutative diagram, where all the arrows are \(T_{\sigma}\)-equivariant :   

\begin{center}
	\begin{tikzpicture}[description/.style={fill=white,inner sep=2pt}]
	\matrix (m) [matrix of math nodes, row sep=3em,
	column sep=1.5em, text height=1.5ex, text depth=0.25ex]
	{  \Gamma(U_{\sigma}, \mathcal{F})  & &  \Gamma(U_{\sigma}, \mathcal{E}) \\
		\mathcal{F}(x_{\sigma}) & &  \mathcal{E}(x_{\sigma})  \\ };
	\path[right hook->]  (m-1-1) edge node[auto] {}(m-1-3);
	\path[right hook->] (m-2-1) edge node[below] {}(m-2-3);
	\path[->] (m-1-1) edge node[auto] {$ev_{x_{\sigma}}$}(m-2-1);
	\path[->] (m-1-3) edge node[auto] {$ev_{x_{\sigma}}$} (m-2-3);
	
	\end{tikzpicture}
\end{center}

Then we can choose a maximal \(T_{\sigma}\)-stable subspace \(F_{\sigma}\) of $\Gamma(U_{\sigma}, \mathcal{F})$ on which $ev_{x_{\sigma}}$ is injective. It can be shown that ${ev_{x_{\sigma}}}|_{F_{\sigma}} : F_{\sigma} \rightarrow \mathcal{F}(x_{\sigma})$ is a \(T_{\sigma}\)-equivariant isomorphism (see \cite[Proposition 2.1.1, STEP 1]{kly}). Since \(F_{\sigma}\) is \(T_{\sigma}\)-stable and ${ev_{x_{\sigma}}}|_{F_{\sigma}}$ is injective, we can choose a maximal \(T_{\sigma}\)-stable subspace \(E_{\sigma}\) of $\Gamma(U_{\sigma}, \mathcal{E})$ containing \(F_{\sigma}\), on which $ev_{x_{\sigma}}$ is injective. Then we have $ev_{x_{\sigma}}|_{E_{\sigma}} : E_{\sigma} \rightarrow \mathcal{E}(x_{\sigma})$ is a \(T_{\sigma}\)-equivariant isomorphism. Moreover we have the following commutative diagram where all the arrows are \(T_{\sigma}\)-equivariant : 

\begin{center}
	\begin{tikzpicture}[description/.style={fill=white,inner sep=2pt}]
	\matrix (m) [matrix of math nodes, row sep=3em,
	column sep=1.5em, text height=1.5ex, text depth=0.25ex]
	{  F_{\sigma}  & &  E_{\sigma} \\
		\mathcal{F}(x_{\sigma}) & &  \mathcal{E}(x_{\sigma})  \\ };
	\path[right hook->]  (m-1-1) edge node[auto] {}(m-1-3);
	\path[right hook->] (m-2-1) edge node[below] {}(m-2-3);
	\path[->] (m-1-1) edge node[auto] {$ev_{x_{\sigma}}$}(m-2-1);
	\path[->] (m-1-1) edge node[left] {$\cong$}(m-2-1);
	\path[->] (m-1-3) edge node[auto] {$ev_{x_{\sigma}}$} (m-2-3);
	\path[->] (m-1-3) edge node[left] {$\cong$} (m-2-3);
	\end{tikzpicture}
\end{center}

Since \(T_{\sigma}\) is reductive and abelian, there exists an eigen basis \(\{s_1, \ldots, s_r\}\) of \(F_{\sigma}\) which extends to an eigen basis \(\{s_1, \ldots, s_r, s_{r+1}, \ldots, s_m\}\) of \(E_{\sigma}\) (see \cite[Theorem 1.23]{brion}). Again these sections are linearly independent at every point \(x \in U_{\sigma}\) (see \cite[Proposition \(2.1.1\), STEP \(2\)]{kly}). Hence we have \(T_{\sigma}\)-equivariant trivializations for the toric vector bundles \(\mathcal{F}\) and \(\mathcal{E}\) respectively denoted by \(\pi_{\sigma}:{\mathcal{E}}|_{U_{\sigma}} \backsimeq  U_{\sigma} \times \mathcal{E}(x_{\sigma})\) and \(\pi^F_{\sigma}:{\mathcal{F}}|_{U_{\sigma}} \backsimeq  U_{\sigma} \times \mathcal{F}(x_{\sigma})\) respectively, where \(\pi_{\sigma}^{-1}((x, \sum_{i=1}^m a_is_i)=\sum_{i=1}^m a_is_i(x) \) for \(x \in U_{\sigma}\) and scalars \(a_1, \ldots, a_m\). Note that the we can give an action of \(T\) on \( U_{\sigma} \times \mathcal{E}(x_{\sigma})\), determined by some representation of the torus \(\phi_{\sigma}: T \rightarrow \text{Aut}(\mathcal{E}(x_{\sigma}))\) extending the action of \(T_{\sigma}\) on $\mathcal{E}(x_{\sigma})$ as follows : the characters of the \(\phi_{\sigma}\) action are specified by taking a preimage of the characters occurring in the \(T_{\sigma}\)-isotypical decomposition of the space \(\mathcal{E}(x_{\sigma})\), under the surjection \(\widehat{T}\rightarrow \widehat{T}_{\sigma} \). We see that, the action \(\phi_{\sigma}\) may not induce the same \(T\)-action on \(E\), which was induced from the \(T\)-action on \(\Gamma(U_{\sigma}, \mathcal{E})\). But by \cite[Proposition 2.2.1(ii)]{kly}),  no matter which extension \(\phi_{\sigma}\) we choose, \({\mathcal{E}}|_{U_{\sigma}}\) is isomorphic to \( U_{\sigma} \times \mathcal{E}(x_{\sigma})\) as \(T\)-equivariant vector bundles. Similarly, the \(T\)-action on \( U_{\sigma} \times \mathcal{F}(x_{\sigma})\) is determined by some extension of the action of \(T_{\sigma}\) on \(\mathcal{F}(x_{\sigma})\) to a representation of the torus \(\phi^F_{\sigma}: T \rightarrow \text{Aut}(\mathcal{F}(x_{\sigma}))\). Further we have  \(\phi^F_{\sigma}(t)={\phi_{\sigma}(t)}|_{\mathcal{F}(x_{\sigma})}\) for all \(t \in T\). Hence we have the following commutative diagram with equivariant arrows :

\begin{center}
	\begin{tikzpicture}[description/.style={fill=white,inner sep=2pt}]
	\matrix (m) [matrix of math nodes, row sep=3em,
	column sep=2.5em, text height=1.5ex, text depth=0.25ex]
	{ {\mathcal{F}}|_{U_{\sigma}}    & &  U_{\sigma} \times \mathcal{F}(x_{\sigma})  \\
		{\mathcal{E}}|_{U_{\sigma}}	 & &    U_{\sigma} \times \mathcal{E}(x_{\sigma})  \\ };
	\path[->]  (m-1-1) edge node[auto] {}(m-1-3);
	\path[->] (m-2-1) edge node[below] {}(m-2-3);
	\path[right hook->] (m-1-1) edge node[auto] {}(m-2-1);
	\path[->] (m-1-1) edge node[below] {$\cong$}(m-1-3);
	\path[->] (m-1-1) edge node[above] {$\pi^F_{\sigma}$}(m-1-3);
	\path[right hook ->] (m-1-3) edge node[auto] {} (m-2-3);
	\path[->] (m-2-1) edge node[below] {$\cong$} (m-2-3);
	\path[->] (m-2-1) edge node[above] {$\pi_{\sigma}$} (m-2-3);
	\end{tikzpicture}   
\end{center}

Now we construct transition functions \(f_{\sigma \tau}: U_{{\sigma} \cap \tau } \rightarrow \text{Hom}(\mathcal{E}(x_{\tau}), \mathcal{E}(x_{\sigma}))\) and \(~ g_{\sigma \tau}: U_{{\sigma} \cap \tau } \rightarrow \text{Hom}(\mathcal{F}(x_{\tau}), \mathcal{F}(x_{\sigma}))\)   for the bundles \(\mathcal{E}\) and \(\mathcal{F}\) respectively, using the above trivializations. Note that from the construction we have \(g_{\sigma \tau}(x)={f_{\sigma \tau}(x)}|_{\mathcal{F}(x_{\sigma})}\) for all \(x \in U_{{\sigma} \cap \tau }\).

Let \(x_0\) denote the identity element of $T$. Let us use the isomorphisms \(f_{\sigma \tau}(x_0):\mathcal{E}(x_{\tau}) \rightarrow  \mathcal{E}(x_{\sigma})\) (resp. \(g_{\sigma \tau}(x_0):\mathcal{F}(x_{\tau}) \rightarrow  \mathcal{F}(x_{\sigma})\) ) to identify all the spaces \(\mathcal{E}(x_{\sigma})\) (resp. \(\mathcal{F}(x_{\sigma})\)) with the fiber \(E=\mathcal{E}(x_{0})\) (resp. \(F=\mathcal{F}(x_{0})\)) in the following way :\\
Define the transition functions as follows :
\[\tilde{g}_{\sigma \tau}:U_{{\sigma} \cap \tau} \rightarrow GL(\mathcal{F}(x_0)) \text{ given by } x \longmapsto \tilde{g}_{\sigma \tau}(x):=g_{o \sigma}(x_0) \circ g_{\sigma \tau}(x) \circ g_{o \tau}(x_0)^{-1} \text{ and}\]
\[\tilde{f}_{\sigma \tau}:U_{{\sigma} \cap \tau} \rightarrow GL(\mathcal{E}(x_0)) \text{ given by } x \longmapsto \tilde{f}_{\sigma \tau}(x):=f_{o \sigma}(x_0) \circ f_{\sigma \tau}(x) \circ f_{o \tau}(x_0)^{-1},\] where \(o\) in the subscript of \(f\) and \(g\) denotes the zero cone in the fan $\Delta$.  

Note that \(\tilde{g}_{\sigma \tau}(x )={\tilde{f}_{\sigma \tau}(x )}|_{\mathcal{F}(x_{0})}\) for all \(x \in U_{{\sigma} \cap \tau} \). Now define the action of \(T\) on \(\mathcal{E}(x_0)\) by \(\tilde{\phi}_{\sigma}(t)=f_{o \sigma}(x_0) \circ \phi_{\sigma}(t) \circ f_{ \sigma o}(x_0)  \) and define an action of \(T\) on \(\mathcal{F}(x_0)\) by \(\tilde{\phi}^F_{\sigma}(t)=g_{o \sigma}(x_0) \circ \phi^F_{\sigma}(t) \circ g_{ \sigma o}(x_0)  \). Note that we have \( \tilde{\phi}^F_{\sigma}(t)={\tilde{\phi}_{\sigma}(t)}|_{F}  \).

Hence for any $\sigma \in \Delta$, $F$ is a $T_{\sigma}$-stable subspace of $E$ where the $T_{\sigma}$-action on $E$ is induced from the map \(\tilde{\phi}_{\sigma}\). Since $T_{\sigma}$ is reductive, the finite dimensional $T_{\sigma}$-module $E$ is semisimple. This implies that the $T_{\sigma}$-stable subspace $F$ of $E$ has a $T_{\sigma}$-stable complement. In other words, there exists an eigen basis of $F$, say $\{f_1, \ldots, f_r \}$ which extends to an eigen basis of $E$, say $B=\{f_1, \ldots, f_r, f'_1, \ldots, f'_m \}$. 
Let \(t \cdot f_i= \chi_\textsubscript{\(i\)}(t)f_i \text{ for } i=1, \ldots, r \text{ and let}\) \(t \cdot f'_j= \chi'_\textsubscript{\(j\)}(t)f'_j\) for \(j=1, \ldots, m \) for all $t \in T_{\sigma}$. Set  $E_{\chi_\textsubscript{\(i\)}}=\text{Span}(f_i)$ and $E_{\chi'_\textsubscript{\(j\)}}=\text{Span}(f'_j)$. Then we get a $T_{\sigma}$-isotypical decomposition \( E=E_{\chi_\textsubscript{\(1\)}} \oplus \cdots \oplus E_{\chi_\textsubscript{\(r\)}} \oplus E_{\chi'_\textsubscript{\(1\)}} \oplus \cdots \oplus E_{\chi'_\textsubscript{\(m\)}}, \) where $\chi_\textsubscript{\(i\)}$ and $\chi'_\textsubscript{\(j\)}$s occurring in the expression need not be distinct. By the compatibility of Klyachko's classification theorem (\cite[Theorem \(2.2.1\)]{kly}) we have $E^{\rho}(i)=\sum\limits_{\langle \psi, v_{\rho} \rangle \geq i} E_{\psi}$ where $\psi$ varies over the set $\{\chi_\textsubscript{\(1\)}, \ldots ,\chi_\textsubscript{\(r\)},  \chi'_\textsubscript{\(1\)}, \ldots ,\chi'_\textsubscript{\(m\)}\}$, and $F^{\rho}(i)=\sum\limits_{\langle \chi, v_{\rho} \rangle \geq i} E_{\chi}$ where $\chi$ varies over the set $\{\chi_\textsubscript{\( 1\)}, \ldots ,\chi_\textsubscript{\(r\)}\}$.

Clearly, $F^{\rho}(i) \subseteq F \cap E^{\rho}(i)$. Let $x \in F \cap E^{\rho}(i)$ and write $x=a_1f_1 + \cdots +a_rf_r= d_1f_1 + \cdots +d_rf_r +d'_1f'_1 + \cdots +d'_mf'_m$, where $d_j \neq 0$ implies $\langle \chi_\textsubscript{\(j\)}, v_{\rho} \rangle \geq i$ and $d'_j \neq 0$ implies $\langle \chi'_\textsubscript{\(j\)}, v_{\rho} \rangle \geq i$. Comparing both the expressions, we get $d_j=a_j$ and $d'_j=0$ which implies that $x \in F^{\rho}(i)$. Hence $F^{\rho}(i)=F \cap E^{\rho}(i)$ holds and the collection of subspaces $\{F, \{E^{\rho}(i) \}_{\rho \in \sigma(1)} \}$ of $E$ forms a distributive lattice.

Conversely, let $F \subseteq E$ and  $F^{\rho}(i)=F \cap E^{\rho}(i)$ be given such that the collection of subspaces $\{F, \{E^{\rho}(i) \}_{\rho \in \sigma(1)} \}$ of $E$ forms a distributive lattice for all $\sigma \in \Delta$. Fixing \(\sigma \in \Delta \), let \(B_{\sigma}=\{f_1, \ldots, f_r, f_{r+1}, \ldots, f_m\}\) be a basis of \(E\) such that \(B_{\sigma} \cap F (=\{f_1, \ldots, f_r\}, \text{say})\) forms a basis of \(F\) and \(B_{\sigma} \cap E^{\rho}(i)\) forms a basis of \(E^{\rho}(i)\) for all \(\rho \in \sigma(1)\). Since the filtrations are decreasing and full, for each \(\rho \in \sigma(1) \) there exists an integer \(n^j_{\rho}\) such that \(f_j \in  E^{\rho}(i)\) for all \(i \leq n^j_{\rho}\) and \(f_j \notin  E^{\rho}(i)\) for all \(i > n^j_{\rho}\). Now define the map \(\chi_\textsubscript{\(j\)} : N_{\sigma} \rightarrow \Z \) by sending \(v_{\rho}\) to \(n^j_{\rho}\) . Since \(\sigma\) is a nonsingular cone \(\chi_\textsubscript{\(j\)}\) becomes a character of \(T_{\sigma}\). Then we can define an action of \(T_{\sigma}\) on \(E\) given by \(t \cdot f_j= \chi_\textsubscript{\(j\)}(t) f_j \) for all \(t \in T_{\sigma}\) and \(j=1, \ldots, m\). Under this action $F$ is stable. Therefore, we get a \(T_{\sigma}\) isotypical decomposition of \(E\) and \(F\) such that the filtrations $\left( F, \{F^{\rho}(i) \}_{\rho \in \Delta(1)} \right)$ and $\left( E, \{E^{\rho}(i) \}_{\rho \in \Delta(1)} \right)$ satisfy the Klyachko's compatibility condition \((\bf{C})\) of \cite[Theorem \(2.2.1\)]{kly}.

From the \(T_{\sigma}\) action on \(E\), we get a map \(T_{\sigma}\rightarrow \text{Aut}(E)\) which can be extended to \(\phi_{\sigma} : T \rightarrow \text{Aut}(E)\). Similarly, we can define \(\psi_{\sigma} : T \rightarrow \text{Aut}(F)\). Note that \(\psi_{\sigma}(t) ={\phi_{\sigma}(t)}|_{F} \). Then define the functions \(f_{\sigma \tau}:T \rightarrow \text{Aut}(E),\) by setting \(f_{\sigma \tau}(t \cdot x_0)=\phi_{\sigma}(t) \phi_{\tau}(t)^{-1}\)for \( \tau, \sigma \in \Delta \). Similarly define \( g_{\sigma \tau }: T \rightarrow \text{Aut}(F)\) by  \(g_{\sigma \tau}(t \cdot x_0)=\psi_{\sigma}(t) \psi_{\tau}(t)^{-1}\) for \( \tau, \sigma \in \Delta \). Since the filtrations satisfy the compatibility condition, these functions extend from \(T\) to \(U_{\sigma} \cap U_{\tau} \) and hence give rise to transition functions which we again denote by \(f_{\sigma \tau}\) and \(g_{\sigma \tau}\) respectively. Note that \({f_{\sigma \tau}(x)}|_{F} =g_{\sigma \tau}(x)\) for all \(x \in U_{\sigma} \cap U_{\tau} \).

The toric vector bundle \(\mathcal{E}\) is given by \( \left(  \amalg_{\sigma \in \Delta} \left(  U_{\sigma} \times E \right) \right)   \diagup \sim  \), where \((x, e) \sim (x,f_{\sigma \tau}(x)e )\) whenever \(x \in U_{\sigma} \cap U_{\tau} \). Similarly, the toric vector bundle \(\mathcal{F}\) is given by \( \left(  \amalg_{\sigma \in \Delta} \left(  U_{\sigma} \times F \right) \right)   \diagup \sim  \), where \((x, e) \sim (x,g_{\sigma \tau}(x)e )\) whenever \(x \in U_{\sigma} \cap U_{\tau} \). For each \(\sigma \in \Delta \) we can define the map \(U_{\sigma} \times F \rightarrow U_{\sigma} \times E\) by sending \((x, e) \longmapsto (x, e)\), where \(x \in U_{\sigma} \) and \(e \in F\). Since \({f_{\sigma \tau}(x)}|_{F} =g_{\sigma \tau}(x)\) for all \(x \in U_{\sigma} \cap U_{\tau} \), the above maps glue together to give a bundle map \( \mathcal{F} \hookrightarrow \mathcal{E} \) showing that \(\mathcal{F}\) is a subbundle of the vector bundle \(\mathcal{E}\). From the equality \(\psi_{\sigma}(t) e =\phi_{\sigma}(t) e \) for all \(t \in T_{\sigma}\) and \(e \in F \), the equivariance of the inclusion map \(\mathcal{F} \hookrightarrow \mathcal{E} \) follows. Thus we conclude that \(\mathcal{F}\) is a toric subbundle of \(\mathcal{E}\).
$\hfill\square$

\begin{cor}
	Let \(X\) be a nonsingular complete toric variety of dimension \(n\). Then the tangent bundle \(T_X\) possesses a toric line subbundle if and only if there exists a ray $\alpha \in \Delta(1)$ such that for any $\sigma \in \Delta(n) $, either \(v_{\alpha}\) or \(-v_{\alpha}\) is a primitive ray generator of $\sigma$.  
\end{cor}
\noindent
{\bf Proof:} Let $(E, \{E^{\rho}(i)\}_{\rho \in \Delta(1)})$ be the associated filtrations  of the tangent bundle \(T_X\) (see Example \ref{tb}). Let $\mathcal{L}$ be a toric line subbundle of \(T_X\) with associated filtrations \((L, \{L^{\rho}(i)\}_{\rho \in \Delta(1)})\) satisfying the condition of Proposition \ref{6}. Let $\sigma=\text{Cone}(v_{\rho_1}, \ldots, v_{\rho_n}) \in \Delta(n)$. Then since the collection of vector spaces \(\{L,E^{\rho_1}(1) =\text{Span}(v_{\rho_1}), \ldots, E^{\rho_n}(1) =\text{Span}(v_{\rho_n} )\}\) forms a distributive lattice, we must have \(L=\text{Span}(v_{\alpha})\) for some $\alpha \in \sigma(1)$. Now for $\sigma' \in \Delta(n)$ other than $\sigma$, again since the collection of vector spaces \(\{L,\{E^{\rho}(1) \}_{\rho \in \sigma'(1)} \}\) forms a distributive lattice, it follows that either \(v_{\alpha}\) or \(-v_{\alpha}\) is a primitive ray generator of $\sigma'$.

Conversely, if the given condition holds, then corresponding to \(L=\text{Span}(v_{\alpha})\) we have a toric subbundle of \(T_X\) by Proposition \ref{6}.
$\hfill\square$

As an application, we see that, the tangent bundle \(T_{\mathbb{P}^n}\) does not have any toric line subbundle. For the Bott tower \(M_k\), corresponding to \(L=\text{Span}(v_{k})\), we have that \(D_k+D_{2k}\) is a toric line subbundle of  \(T_{M_k}\).

\begin{cor}\label{14}
Consider  the tangent bundle $T_X$ of the Hirzebruch surface $X=\mathcal{H}_r$ with the polarization $H=aD_3+bD_4$ $(a, b >0)$. Then for $r \geq 2$, $T_X$ is unstable.
\end{cor}
\noindent
{\bf Proof:} Note that \(c_1(T_X)=D_1+D_2+D_3+D_4\). We have the following relations $D_1 \cdot D_3=D_3 \cdot D_3=D_2 \cdot D_4=0$, $D_2 \cdot D_3=D_4 \cdot D_3=D_1 \cdot D_4=1$ and $D_4 \cdot D_4=r$.
Hence the slope is given by
 $$\mu (T_X)= \frac{\text{deg }(T_X)}{\text{rank }(T_X)}= \frac{2a+b(r+2)}{2}.$$ 
From the above discussion, we see that $D_2+D_4$ is a toric line subbundle of $T_X$ with slope \(2a +br\).  Thus $\mu(D_2+D_4 ) \leq \mu(T_X)$ if and only if $rb \leq 2b-2a$. Hence for $r \geq 2$, $T_X$ is unstable. 
$\hfill\square$

\bibliographystyle{plain}
\bibliography{refs_paper}

\end{document}